\documentclass{amsart}

\newtheorem{theorem}{Theorem}[section]
\newtheorem*{main result}{Theorem}
\newtheorem*{claim}{Claim}
\newtheorem{lemma}[theorem]{Lemma}
\newtheorem{proposition}[theorem]{Proposition}

\theoremstyle{definition}
\newtheorem{definition}[theorem]{Definition}
\newtheorem{example}[theorem]{Example}
\newtheorem{remark}[theorem]{Remark}
\newtheorem{remarks}[theorem]{Remarks}
\newtheorem*{acknowledgement}{Acknowledgement}

\usepackage{amscd,amssymb,epsfig,amsfonts}

\newcommand{\NN}{\mathbb{N}}
\newcommand{\ZZ}{\mathbb{Z}}
\newcommand{\QQ}{\mathbb{Q}}
\newcommand{\RR}{\mathbb{R}}

\newcommand{\PP}{\mathbb{P}}
\renewcommand{\AA}{\mathbb{A}}
\newcommand{\GG}{\mathbb{G}}

\newcommand  {\shB}     {\mathcal{B}}
\newcommand  {\shC}     {\mathcal{C}}

\newcommand  {\shExt}   {\mathcal{E} \!\text{\textit{xt}}}
\newcommand  {\shE}     {\mathcal{E}}
\newcommand  {\shF}     {\mathcal{F}}
\newcommand  {\shG}     {\mathcal{G}}

\newcommand  {\shHom}   {\mathcal{H}\!\text{\textit{om}}}

\newcommand  {\shM}     {\mathcal{M}}

\newcommand  {\shN}     {\mathcal{N}}
\newcommand  {\shL}     {\mathcal{L}}

\newcommand  {\shS}     {\mathcal{S}}
\newcommand  {\shT}     {\mathcal{T}}

\newcommand  {\foX}     {\mathfrak{X}}

\newcommand  {\Art}     {\operatorname{Art}}

\newcommand  {\Aut}     {\operatorname{Aut}}

\newcommand  {\Div}     {\operatorname{Div}}
\newcommand  {\dlog}    {\operatorname{dlog}}

\newcommand  {\Ext}     {\operatorname{Ext}}

\newcommand  {\Fr}      {\operatorname{Fr}}

\newcommand  {\Hom}     {\operatorname{Hom}}
\newcommand  {\Hilb}     {\operatorname{Hilb}}

\newcommand  {\id}      {\operatorname{id}}

\renewcommand  {\k}     {\kappa}
\renewcommand  {\ker }  {\operatorname{kern}}
\newcommand  {\coker }  {\operatorname{cokern}}

\newcommand  {\LD}      {\operatorname{LD}}

\newcommand  {\lra}     {\longrightarrow}

\newcommand  {\maxid}   {\mathfrak{m}}

\newcommand  {\N}       {\operatorname{N}}

\newcommand  {\NS}      {\operatorname{NS}}

\renewcommand{\O}       {\mathcal{O}}

\newcommand  {\Pic}     {\operatorname{Pic}}

\newcommand  {\quadand} {\quad\text{and}\quad}

\newcommand  {\ra}      {\rightarrow}

\newcommand  {\Set}     {\operatorname{Set}}

\newcommand  {\Sing}    {\operatorname{Sing}}

\newcommand  {\Spec}    {\operatorname{Spec}}
\newcommand  {\Spf}     {\operatorname{Spf}}

\def\mydate{\number\day\space\ifcase\month \or January\or February\or 
March\or April\or May\or\June\or\July\or
August\or September\or October\or November\or December\fi \space\number\year}

\begin{document}

\title[Enriques surfaces in characteristic two]
{Logarithmic deformations of normal crossing 
Enriques surfaces in characteristic two}
\author[Stefan Schroeer]{Stefan Schr\"oer}
\address{Mathematische Fakult\"at, Ruhr-Universit\"at, 
               44780 Bochum, Germany}

\curraddr{M.I.T. Department of Mathematics, 
          77 Massachusetts Avenue, Cambridge MA 02139-4307, USA}

\email{s.schroeer@ruhr-uni-bochum.de}



\dedicatory{25 March 2001}

\begin{abstract}
Working in characteristic two, I classify 
nonsmooth Enriques surfaces with  normal crossing singularities.
Using Kato's theory  of logarithmic structures, 
I show that   such surfaces are smoothable and lift to characteristic zero,  
provided they are 
d-semistable. 
\end{abstract}

\maketitle

\section*{Introduction}

Smooth Enriques surfaces occupy a special position in the classification 
of surfaces:
They are   similar to rational surfaces but have Kodaira dimension zero.
Here we shall study   Enriques surfaces that are normal crossing 
rather than smooth:
What they are; how they look like; their fibrations or embeddings;
and whether they deform to smooth surfaces or lift to characteristic zero.
Since $p=2$ 
 is the most exciting prime with respect to Enriques surfaces and to surface 
normal crossings as well,
I shall restrict my attention to characteristic 
two.

Little seems to be known on the relation of Enriques surfaces in 
characteristic two and  characteristic zero. However, Ekedahl and 
Shepherd--Barron recently announced
that smooth Enriques surface in characteristic two lift to 
characteristic zero.  Our main results is
the following Theorem: 

\begin{main result}
Nonsmooth   d-semistable Enriques surfaces with normal crossings
are smoothable and lift to characteristic zero.
\end{main result}

Note that this does not require projectivity. The 
result might be useful for the construction of
complete moduli spaces over $\Spec(\ZZ)$ for Enriques surfaces.

Working over the complex numbers, 
Kulikov \cite{Kulikov 1977}  started the study of simple normal crossing 
K3 and Enriques surfaces.
Using Hodge theory, he obtained a classification of such surfaces. 
Taking the peculiarities of characteristic 
$p=2$ into account, I shall give a similar classification.  Criteria for   
classification are:
Structure of Picard scheme (classical, ordinary, or supersingular);  
multiplicity   of singularities (type II  or type III); and   nature of 
irreducible components
(simple or nonsimple).

Friedman's paper \cite{Friedman 1983} is fundamental for
the theory of smoothings. Working in the complex analytic category,
he showed that precisely the  d-semistable K3 surfaces from Kulikov's list 
deform to smooth K3 surfaces.
Using logarithmic structures, Kawamata and Namikawa 
\cite{Kawamata; Namikawa 1994}  
simplified Friedman's   proof. Little attention,  however, was paid to positive
characteristics.

The theory of logarithmic structures goes back to Fontaine and Illusie
and was developed by K.~Kato \cite{Kato 1989}. 
I shall apply the work of 
F.~Kato \cite{Kato 1996} on log deformations  to d-semistable Enriques
surfaces. We shall see that the   log deformation functor
is formally smooth for classical and ordinary Enriques surfaces. In contrast, 
supersingular surfaces are obstructed.

In the   category of schemes  the existence of   formal smoothings
does not  imply the existence of algebraic smoothings.
In light of Grothendieck's Existence Theorem, algebraization hinges on the
presence of an   formal ample sheaf. 
There are no obstructions for for classical Enriques surfaces. 
For type III Enriques surfaces,
however, we first have to pass to minus-one-form  and check that a 
sufficiently large part
of the versal deformation is algebraizable.

\begin{acknowledgement}
I wish to thank Bernd Siebert and Hubert Flenner for   stimulating discussions.
\end{acknowledgement}

\section{Enriques surfaces with normal crossings}
\label{Enriques surfaces with normal crossings}

Fix a ground field 
$k$, for the moment of arbitrary characteristic 
$p\geq 0$.
To study degenerations of smooth Enriques surfaces we have to  decide on
a  suitable notion of singular Enriques surfaces.

\begin{definition}
\label{normal crossing Enriques surface}
A   proper 
$k$-surface $X$ with 
$k=\Gamma(X,\O_X)$ is called   a \emph{normal crossing Enriques surface}
if  the singularities of 
$X$ are normal crossings, and the canonical class
$K_X $ is numerically trivial, and $\chi(\O_{X})=1$ holds.
\end{definition}

For smooth Enriques surfaces, this coincides with the usual  
definition $K_X\equiv 0$ and $b_2(X)=10$
(compare \cite{Cossec; Dolgachev 1989}, pp.\ 72--74).
There are interesting degenerations of Enriques
surfaces with $K_X\not\equiv 0$ 
(so-called flower pot degenerations). For this paper, however, I am
content with the case $K_X$ numerically trivial.

Note that  normal crossing is a local condition, 
such that the irreducible components of $X$ might be
nonsmooth. Call $X$
\emph{simple} if all its irreducible components are smooth, 
and \emph{nonsimple} otherwise.

By definition, the singularities  are normal crossings if for each point 
$x\in X$ there is an isomorphism
$\O_{X,x}^\wedge\simeq \kappa(x)[[T_1,T_2,T_3]] /( T_1\ldots T_l)$. 
The integer 
$l\in\left\{  1,2,3 \right\}$ is called the \emph{multiplicity} of 
$x\in X$.
A point  $x\in X$ with multiplicity $l=1$ is   smooth.  Points
with multiplicity $l=2$ or $l=3$ are called \emph{double points} or
\emph{triple points}, respectively.
We shall say that $X$ is  an Enriques surface of \emph{type II} if 
it is a normal crossing Enriques surface
with double points but without triple points; and of \emph{type III} 
if it contains triple points.

Another way to distinguish normal crossing Enriques surfaces 
involves the Picard scheme.
Let 
$\Pic^\tau(X)$ be the group of numerically trivial line bundles and 
$\Pic^\tau_X$ be the corresponding group scheme. 
These objects behave like in the smooth case:

\begin{proposition}
\label{Picard group}
Suppose 
$X$ is a normal crossing Enriques surface. Then 
$\Pic^\tau_X$ is an affine group scheme of length 2. Furthermore,
$\Pic^\tau(X)$ is generated by 
$K_X$.
\end{proposition}

\proof
Note that a numerically trivial invertible $\O_X$-module is trivial 
if and only if
it has a nonzero section, because $X$ is connected and reduced. Suppose  
$\shL\neq\O_{X}$ is such a line bundle. By Riemann--Roch,  
$\chi(\shL)=\chi(\O_{X})=1$. On the other hand, 
$$
\chi(\shL) \leq h^0( \shL) + h^0(\shL^\vee\otimes\omega_X) = 
h^0(\shL^\vee\otimes\omega_X),
$$
consequently  $\shL\simeq\omega_X$. Hence, if 
$\Pic^\tau(X)$ is nontrivial, it    is generated by 
$K_X$ and has order two.

In light of this, the connected component 
$\Pic^0_X\subset\Pic^\tau_X$ is discrete. Its Lie algebra is
$H^1(X,\O_{X})$. There are two possibilities: 
If 
$K_X\neq 0$, then
$h^1(\O_{X})=0$, so   
$\Pic^0_X=0$ and 
$\Pic^\tau_X=\ZZ/2\ZZ$.
If
$K_X=0$, then 
$\Pic^{0}_X=\Pic^\tau_X$. Moreover, $h^1(\O_X)=1$, so $\Pic^0_X$ has length
two. In both cases we see that 
$\Pic^\tau_X$ is an affine group scheme of length two.
\qed

\medskip
Besides the \'etale group scheme $\ZZ/2\ZZ=\Spec[T\mid T^2=T]$, 
there are two radical (that is, purely inseparable)  group schemes of 
length two in
characteristic 
$p=2$:
$$
\mu_2 = \Spec k[T\mid T^2=1] \quadand
\alpha_2 =\Spec k[T\mid T^2 =0].
$$
Note that 
$\alpha_2$ and $\ZZ/2\ZZ $ are unipotent group schemes, whereas 
$\mu_2 $ is a multiplicative group scheme.
We shall say that
$X$ is a \emph{classical} Enriques surface if $X$ is a 
normal crossing Enriques surface with 
$\Pic^\tau_X=\ZZ/2\ZZ$;  \emph{ordinary} if 
$\Pic^\tau_X=\mu_2$; and \emph{supersingular} if 
$\Pic^\tau_X=\alpha_2$.

\section{How do normal crossing Enriques surfaces look like?}
\label{How do normal crossing Enriques surfaces look like?}

Fix an algebraically closed ground field 
$k$ of characteristic
$p=2$, and let 
$X$ be an Enriques surface of  type II or III.
The task now is to determine the structure of such surfaces. The idea is
to reconstruct $X$ from its normalization via gluing.

Let 
$\nu:S\ra X$ be the normalization, 
$C\subset S$ the reduced ramification locus of the  normalization map,
and 
$D\subset X$ the reduced singular locus. Then 
$S$ is a smooth surface, 
$C\subset S$ is a normal crossing divisor, and $D$ is a seminormal curve. 
The commutative diagram 
$$
\begin{CD}
C     @>>> S    \\
@V\varphi VV     @VV\nu V\\
D    @>>>  X 
\end{CD}
$$
is cartesian and cocartesian.  Hence we can recover 
$X$ from the smooth surface 
$S$ and the \emph{gluing map} 
$\varphi:C\ra D$. The diagram yields a short exact sequence
$$
1 \lra \O_{X}^\times  \lra \O_{S}^\times\oplus \O_{D}^\times  
\lra \O_{C}^\times  \lra 1,
$$
which in turn gives a long exact sequence of abelian sheaves
\begin{equation}
\label{Picard sequence}
H^0(\O_C^\times) \lra \Pic(X)   \lra \Pic(S)\oplus\Pic(D)  \lra \Pic(C) 
\lra H^2(\O_{X}^\times).  
\end{equation}
The relative dualizing sheaf 
$\omega_{S/X}$ coincides with the conductor ideal of the inclusion 
$\O_X\subset \nu_*(\O_S)$, so we have 
$$
K_S= K_{S/X}+\nu^*(K_X)\equiv -C \quadand 2K_S=-2C.
$$
Decompose 
$S=S_1\cup\ldots\cup S_n$ into connected components, and let 
$X=X_1\cup\ldots\cup X_n$ be the corresponding irreducible components.
Then
$2K_{S_i}= -2C_{S_i}$. By the   classification of surfaces \cite{Mumford
1969}, each component 
$S_i$ is ruled.

The ramification locus 
$C\subset S$ is a divisor with normal crossings.  Its global structure
is easy: The adjunction formula gives 
$2K_{C}= 0$. It follows $K_C=0$, so   
$C$ is a disjoint union of elliptic curves and cycles of rational curves.
Here a \emph{cycle of rational curves}  means a seminormal curve
isomorphic to   
$\PP^1\times\ZZ/m\ZZ$ (with 
$m\geq 1$) modulo the relation  
$(\infty,i)\sim(0,i+1)$ for $i\in\ZZ/m\ZZ$.

\begin{proposition}
\label{ramification on rational component}
Suppose the component
$S_i\subset S$ is  a rational surface. Then the ramification locus
$ S_i\cap C$ is   either an elliptic curve  or a cycle of
rational curves.
\end{proposition}

\proof
Set 
$C_i=S_i\cap C$.
Then 
$K_{S_i}=-C_i$ because 
$\Pic(S_i)$ is torsion free.
The exact sequence 
$$
0 \lra H^0(S_i,\O_{S_i})  \lra H^0(C_i,\O_{C_i} ) \lra  H^1(S_i,\omega_{S_i})
$$
and 
$H^1(S_i,\O_{S_i})=0$ implies that 
$C_i$ is connected, so $C_i$ is as desired.
\qed

\medskip
To proceed, we need a fact on elliptic curves  in characteristic $p=2$.

\begin{lemma}
\label{Hasse}
For an elliptic curve $E$ over $k$, the following  are equivalent.
\renewcommand{\labelenumi}{(\roman{enumi})}
\begin{enumerate}
\item The Frobenius map 
$\Fr^*:H^1(E,\O_E)\ra H^1(E,\O_E)$ is zero.
\item 
There is an inclusion of group schemes
$\alpha_2\subset E$. 
\item
$E$ has a Weierstrass equation of the form $y^2 + y = x^3$.
\item
The $j$-invariant is $j(E)=0$.
\item
The group 
$E(k)$ has no 2-torsion.
\item  $\Aut(E)$ is a group of  order 24.
\end{enumerate}
\end{lemma}

For a proof, see Silverman \cite{Silverman 1986}, Appendix A. 
Note that each  elliptic curve over $k$ has a Weierstrass equation 
$y^2 + \alpha xy + y = x^3$ with $\alpha^3\neq 1$, which
is  called  
\emph{Deuring normal form}.
The corresponding 
$j$-invariant is 
$j=\alpha^{12}/(\alpha^3-1)$.
Condition (i) means that $E$ has \emph{Hasse invariant} zero. 
Such elliptic curves
are   called \emph{supersingular}. We shall see that supersingular 
elliptic curves are   closely
related to supersingular Enriques surfaces.

\begin{proposition}
\label{ramification on elliptic components}
Suppose the component
$S_i\subset S$ is a  nonrational surface. Then 
there is a unique ruling 
$f_i:S_i\ra B_i$ over an elliptic curve 
$B_i$. 
The   curve 
$C\cap S_i$ is either the   union    of two disjoint
sections, or an elliptic curve double-covering 
$B_i$.
Moreover,  
$\nu^*(K_X)|_{S_i}\neq 0$   if and only if  
$C\cap S_i\ra B_i$ is radical and $j(B_i)\neq 0$.
\end{proposition}

\proof
Set 
$C_i=C\cap S_i$.
The Albanese morphism yields a unique ruling 
$f_i:S_i\ra B_i$ over a smooth nonrational curve 
$B_i$. The Hurwitz Formula applied to 
$C_i\ra B_i$ tells us that 
$B_i$ is elliptic. By L\"uroth's Theorem, 
$C_i$ must be a disjoint union of elliptic curves.
By the adjunction formula,  the projection 
$C_i\ra B_i$ has degree 2, so   
$C_i$ is either the   union of two disjoint sections or
an elliptic curve double-covering 
$B_i$.

It remains to verify the assertion concerning
the numerically trivial invertible 
$\O_{S_i}$-module $\shL_i=\nu^*(\omega_X)|_{S_i}$. The exact
sequence 
$$
0 \lra \omega_{S_i}\otimes \shL_i^\vee  \lra \O_{S_i}  \lra \O_{C_i}  \lra 0
$$
gives an exact sequence 
\begin{equation}
\label{cohomology sequence}
H^1(S_i,\O_{S_i}) \lra H^1(C_i,\O_{C_i})  \lra 
H^2(S_i, \omega_{S_i}\otimes \shL_i^\vee)\lra 0.
\end{equation}
Note that
$H^1(\O_{B_i})\ra H^1(\O_{S_i})$ is bijective, and 
$H^2( \omega_{S_i}\otimes \shL_i^\vee) \simeq H^0(\shL_i)$.
If
$C_i$ is the union of two disjoint section, then
$h^1(\O_{C_i})=2$, hence   
$\shL_i$ has a nontrivial section, so
$\nu^*(K_X)|_{S_i}=0$.

Now suppose that 
$C_i$ is connected.
Choosing  suitable group structures, we may assume that 
$C_i\ra B_i$ is a homomorphism of group  schemes. We obtain an exact
sequence 
$$
0 \lra G  \lra C_i  \lra B_i  \lra 0,
$$
where the kernel 
$G$ is one of  
$ \mu_2, \alpha_2$ or 
$\ZZ/2\ZZ$. Applying the functor $\Ext^*(\cdot,\GG_m)$ to the 
preceding exact sequence,
we obtain another  exact sequence 
$$
0 \lra D(G)  \lra \Pic^0_{B_i}  \lra \Pic^0_{C_i}  \lra 0.
$$
Here 
$D(G)=\Hom(G,\GG_m)$ is the \emph{Cartier dual} of 
$G$. We shall distinguish two cases:
First, suppose 
$G$ is either 
$\ZZ/2\ZZ$ or 
$\alpha_2$. Then  the Cartier dual $D(G) $ is infinitesimal, so the map 
$H^1(\O_{B_i})\ra H^1(\O_{C_i})$ is zero, hence 
$\shL_i^\vee$ has a nontrivial section.  Using the exact sequence  
(\ref{cohomology sequence}) we conclude  
$\nu^*(K_X)|_{S_i}=0$.
Second, suppose 
$G=\mu_2$ is multiplicative.  Then 
$D(G)=\ZZ/2\ZZ$ is \'etale,  so the map 
$H^1(\O_{B_i})\ra H^1(\O_{C_i})$ is bijective. We deduce   
$\nu^*(K_X)|_{S_i}\neq 0$. 

Finally, note that,  if 
$C_i\ra B_i$ is radical (that is, purely inseparable), the curves 
$C_i$ and 
$B_i$ are isomorphic as   
$\ZZ$-schemes without 
$k$-structures, and the projection 
$C_i\ra B_i$ is nothing but the Frobenius morphism 
$\Fr:B_i\ra B_i$. Using Lemma \ref{Hasse}, we see that
$j(B_i)\neq 0$ holds if and only if $G=\mu_2$.
\qed

\medskip
We shall use a \emph{bicoloured graph}  to describe the combinatorics of 
type II surfaces. The
 vertices and edges of this graph  are as follows:  
\begin{itemize}
\item White vertices: $\pi_0(S)$.
\item Black vertices: $\pi_0(D)$.
\item Edges: $\pi_0(C)$.
\end{itemize}
Here $\pi_0(.)$  denotes the set of connected components.
An edge $C'\subset C$ connects a white vertex $S'\subset S$ with a
black vertex $D'\subset D$ if and only if $C'\subset S'$ and 
$\varphi(C')\subset D'$.

\begin{remark}
\label{bicoloured graph}
The preceding definition makes sense for any seminormal  scheme $X$,
with $C\subset S$ the ramification locus of the normalization   $\nu:S\ra X$, 
and $D\subset X$ support of
$\nu_*(\O_S)/\O_X$. For seminormal curves, the construction is due to
Deligne and Rapoport (\cite{Deligne; Rapoport 1972}, section 3.5).
\end{remark}

\begin{theorem}
\label{nonsimple type II}
Suppose 
$X$ is a nonsimple   Enriques surface of type II.
Then the   bicoloured graph 
$\Gamma(X)$ has the form:\\[2ex]
\mbox{}\hfill\epsffile{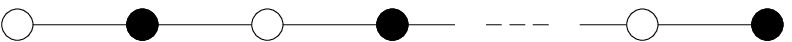}\hfill\mbox{}\\[1ex]
\mbox{}\hspace{6.2em} $S_1$ \hspace{5.4em} $S_2$ \hspace{8.8em} $S_n$\\[1ex]
The first component 
$S_1$ is a rational surface, whereas the  subsequent components
$S_2,\ldots, S_n$ are elliptic ruled.
Furthermore, 
$X$ is ordinary.
\end{theorem} 
 
\proof
The number of rational components 
$S_i\subset S$ equals  
$\chi(\O_S)$. Since 
$X$ is triple-point-free, both the ramification locus $C$ 
and the singular locus 
$D$ are disjoint unions of elliptic curves.
Hence we have 
$$
1=\chi(\O_X)= \chi(\O_S) +\chi(\O_D)-\chi(\O_C) = \chi(\O_S).
$$
Consequently, there is precisely one rational component, say 
$S_1\subset S$. It has at most one neighbor, because 
$C_1$ is connected by Proposition 
\ref{ramification on rational component}. For 
$2\leq i\leq n$, the curve 
$C_i$ has at most two connected components, so the graph
$\Gamma(X)$ is a chain. 
By assumption, 
$X$ has a nonsmooth component, hence the last vertex is black.

It remains to check that 
$X$ is ordinary. Let 
$D'\subset D$ be the connected component corresponding to the rightmost
black vertex, and let $C'\subset C$ be its preimage.
Choosing suitable base points, we may assume that the gluing is given
by an exact sequence
$$
0\lra \ZZ/2\ZZ \lra C'\lra D'\lra 0
$$
of elliptic curves with group structure. The functor $\Ext^*(\cdot,\GG_m)$ 
gives a
dual exact sequence
$$
0\lra D(\ZZ/2\ZZ) \lra \Pic^0_{D'}\lra \Pic^0_{C'}\lra 0.
$$
Using the long exact sequence in (\ref{Picard sequence}),
you easily infer $\mu_2\subset\Pic^\tau_X$. Consequently $X$ is
ordinary. 
\qed

\begin{example}
\label{irreducible}
The case $n=1$ is allowed: Take $S=\PP^2$, and let 
$C\subset S$ be a smooth cubic. The Jacobian 
$\Pic^0_C$ acts freely on 
$C$. Choose an invertible 
$\O_C$-module of order two, and let $\iota:C\ra C$ be the corresponding
free involution, say with quotient $D=C/\iota$. Now $X= S\coprod_C D$  is an 
irreducible
Enriques surface of type II.
\end{example}

\begin{theorem}
\label{simple type II}
Suppose 
$X$ is a simple    Enriques surface of type II.
Then the bicoloured graph 
$\Gamma(X)$ has the form:\\[2ex]
\mbox{}\hfill\epsffile{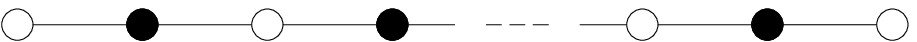}\hfill\mbox{}\\[1ex]
\mbox{}\hspace{4.2em} $S_1$ \hspace{5.5em} $S_2$ 
\hspace{8.8em} $S_{n-1}$\hspace{4.8em} $S_n$\\[1ex]
The first component 
$S_1$ is rational, whereas  
$S_2,\ldots, S_n$ are elliptic ruled.
Moreover, if $f_n:S_n\ra B_n$ is the   ruling on the last component and
$C_n\subset C$ is the component corresponding to the last edge, we have
the following:
\renewcommand{\labelenumi}{(\roman{enumi})}
\begin{enumerate}
\item
The surface $X$ is ordinary if and only if 
$f_n:C_n\ra B_n$ is \'etale.
\item
$X$ is classical if and only if the projection
$C_n\ra B_n$ is radical and 
$j(B_n)\neq 0$.
\item
$X$ is supersingular if and only if $C_n\ra B_n$ is radical and 
$j(B_n)=0$.
\end{enumerate}
\end{theorem}

\proof
You verify the first two assertions as in the preceding proof.
For the last statement, consider for example the case that
$C_n\ra B_n$ is radical with
$j(B_n)=0$. Then $C_n$ and $B_n$ are isomorphic as
$\ZZ$-schemes without $k$-structures, so we easily obtain an exact sequence
$$
0\lra \alpha_2 \lra \Pic^0_{S_n}\lra \Pic^0_{C_n}\lra 0.
$$
Using the exact sequence in (\ref{Picard sequence}),
we infer $\alpha_2\subset \Pic^\tau_X$, so $X$ is supersingular.
You handle the remaining cases in  a similar way.
\qed

\medskip
The next task is to treat type III surfaces. We shall describe their
combinatorics   in terms of  \emph{cell decompositions} $\Pi(X)$ of
compact real 2-manifolds as follows: 
The 1-skeleton of the cell decomposition is the 
bicoloured graph 
$\Gamma(D)$ attached to the seminormal curve $D$ as in Remark
\ref{bicoloured graph}. 
The 2-cells  correspond to  
$\pi_0(S)$. It remains  to specify attaching maps. 
We shall  see that each $C_i=S_i\cap C$ is a cycle  of rational curves, so  
$\Gamma(C_i)$ triangulates the real circle 
$\RR/\ZZ$. For each   2-cell corresponding to a connected component 
$S_i\subset S$, 
identify its boundary circle with
$\Gamma(C_i)$ and use the canonical map $\Gamma(C_i)\ra\Gamma(D)$ for
attaching the 2-cell  to the 1-skeleton $\Gamma(D)$. 
You directly check that this yields a cell
decomposition $\Pi(X)$ of a compact real 2-manifold.

\begin{remarks}
(i) The cell decomposition $\Pi(X)$ is nothing but a \emph{dessin d'enfant}, 
except that we do not require orientability. Such
objects were introduced by Grothendieck \cite{Grothendieck 1997} in his 
anabelian study of the
absolute  Galois group 
$\Aut(\bar{\QQ})$.

(ii) The cell decomposition is uniquely determined  by a 
regular neighborhood of the  1-skeleton 
$\Gamma(D)\subset\Pi(X)$, which 
is a \emph{ribbon graph}. You can specify    ribbon graphs in the following
 way: Choose an  immersion
of 
$\Gamma(D)$ into the real plane and mark those edges whose twisted 
ribbons shall have a single twist.
For an illustration, see Example \ref{example type III}.
\end{remarks}

Now we are ready to describe  surfaces with triple points. The following 
result was
obtained by Kulikov \cite{Kulikov 1977} for simple
normal crossings in the complex analytic case   
via Hodge theoretic arguments.

\begin{theorem}
\label{type III}
Suppose 
$X$ is a type III   Enriques surface. Then the real 2-manifold underlying  
$\Pi(X)$ is the Klein bottle
$\RR\PP^2$.
Each component 
$S_i\subset S$  is rational, and each ramification curve 
$C\cap S_i$ is a cycle of rational curves.
Furthermore, 
$X$ is ordinary.
\end{theorem}

\proof
Since 
$X$ contains a triple point, there is at least on cycle of rational curves 
inside 
$C$. The  component 
$S_j$ containing it is rational, and   
$C\cap S_j$ is a cycle of rational curves.
Consequently, each component 
$S_i$ intersecting 
$S_j$ is also   rational and 
$C_i$ is a cycle of rational curves.
Since 
$X$ is connected, this holds for all components 
$S_i$.
The exact sequence in (\ref{Picard sequence})
yields an inclusion 
$\Pic^\tau_X\subset \Pic^0_D$. Since 
$D$ is seminormal, 
$\Pic^0_D$ is a torus, so the subgroup 
$\Pic^\tau_X$ is multiplicative. Since 
$\ZZ/2\ZZ$ and 
$\alpha_2$ are not multiplicative  in  characteristic two, we infer
that 
$\Pic^\tau_X=\mu_2$, so  
$X$ is ordinary.

It remains to determine the real 2-manifold underlying
$\Pi(X)$. We do this by calculating the Euler characteristic
$\chi=v-e+f$.  
The number of faces 
$f$ equals the number 
$n$ of components 
$S_i$, which also coincides with 
$h^1(\O_C)$.
The difference of vertices and  edges is
$$
v-e =\chi(\Gamma(D)) = h^0(\Gamma(D),k)-h^1(\Gamma(D),k) = 1-h^1(\O_D).
$$
Here we use $h^1(\Gamma(D),k) =  h^1(\O_D)$, which follows from the
exact sequence 
$$
 H^0(\tilde{D},\O_{\tilde{D}}) \oplus H^0(D',\O_{D'}) \lra
H^0(\nu^{-1}(D'),\O_{\nu^{-1}(D')}) \lra H^1(D,\O_{D})\lra 0
$$
for the normalization $\nu:\tilde{D}\ra D$, with $D'=\Sing(D)$.
The exact sequence 
$$
0 \lra H^1(X,\O_{X})  \lra H^1(D,\O_{D})   \lra H^1(C,\O_{C})   
\lra H^2(X,\O_{X}) \lra  0
$$
gives 
$h^1(\O_D)=h^1(\O_C)$. The upshot of this is
$$
\chi = (v-e) + f = 1-h^1(\O_D) + h^1(\O_C) = 1.
$$
By the classification of compact real 2-manifolds,  
$\Pi(X)$ is a cell decomposition of   the Klein bottle
$\RR\PP^2$.
\qed

\begin{example}\label{example type III}
Let me discuss the  easiest type III surfaces, which is
 nonsimple   with
$X=X_1\cup X_2$   two irreducible components. The
following picture illustrates the gluing map $\varphi:C\ra D$:\\[1.5ex]
\mbox{}\hfill\epsffile{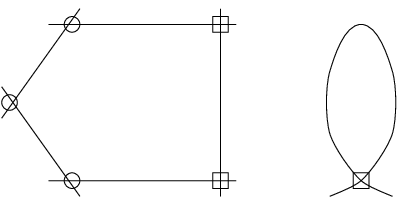}\hfill\epsffile{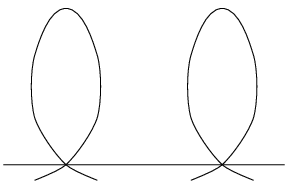}\hfill\mbox{}\\[1ex]
\mbox{}\hspace{10em} $C=C_1\cup C_2$ \hspace{8.7em} $D\subset X$ \\[1.5ex]
The symbols on the intersection points indicate the gluing. The component
$S_1\subset S$ could be the blowing-up of a Hirzebruch surface
$\bar{S}_1$
with an anticanonical 4-cycle of rational curves 
$\bar{C}_1\subset \bar{S}_1$ (two disjoint sections and two fibers),
so that the  center  of the  blowing-up is a singularity of 
$\bar{C}_1$. For 
$C_2\subset S_2$, choose a nodal cubic in 
$\PP^2$. 
The cell decomposition $\Gamma(X)$ is given by an immersion of
the 1-skeleton $\Gamma(D)$:\\[1.5ex]
\mbox{}\hfill\epsffile{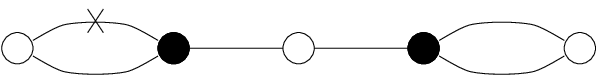}\hfill\mbox{} \\[1.5ex]
Here the unmarked edges  have   untwisted ribbons, whereas the marked edge
has a twisted ribbon.
There is another Enriques surfaces with two irreducible components, such that 
$C_1$ is a 4-cycle, and 
$C_2$ is a 2-cycle of rational curves.
\end{example}

\begin{remark}
As in the smooth case, there is a  canonical double covering for a 
normal crossing Enriques surfaces 
$X$: According to Raynaud \cite{Raynaud 1970}, Proposition 6.2.1, 
the inclusion 
$\Pic^\tau_X\subset\Pic_X$ corresponds to a nontrivial 
principal homogeneous 
$G$-space
$r:\tilde{X}\ra X$, where 
$G=\Hom(\Pic^\tau_X, \GG_m)$ is the Cartier dual.  The surface 
$\tilde{X}$ is locally of complete intersection and has 
$\Gamma(\tilde{X},\O_{\tilde{X}})=k$. Moreover, 
$H^1(\tilde{X},\O_{\tilde{X}})=0$ and 
$K_{\tilde{X}}=0$. This suggests to  call 
$\tilde{X}$ the \emph{K3-like covering} of 
$X$. 

For ordinary Enriques surfaces, 
$G=\ZZ/2\ZZ$ is \'etale, so the K3-like covering
$\tilde{X}$ is a normal crossing  K3 surface.
Suppose that  
$X$ is classical or supersingular. Then 
$G$ is radical, so the covering
$r:\tilde{X}\ra X$ is radical. With the notation of Theorem 
\ref{simple type II}, 
you can easily see that the induced homogeneous
$G$-space 
$r_i:\tilde{X_i}\ra X_i$ is trivial for  
$i=1,\ldots n-1$. Note that the K3-like covering $\tilde{X}$ 
is nonreduced. 
\end{remark}

\section{Projectivity  and d-semistability}

Over the complex numbers, each smooth analytic Enriques surface is projective 
(\cite{Barth; Peters; Van de Ven 1984}, p.\ 184). 
Are  normal crossing Enriques surface projective? In
this section, we shall analyze this problem for type III surfaces $X$.
Let 
$C_j\subset C$ be the irreducible components of the ramification curve 
$C\subset S$ on the normalization $S$. 
Following Miranda and Morrison \cite{Miranda; Morrison 1983}, we say that 
$X$ is in\emph{ minus-one-form} if 
$$
C_j^2 = \left\{
\begin{array}{l}
-1 \quad \text{if $C_j\subset S$ is smooth,}\\
+1 \quad \text{if $C_j\subset S$ is nodal,}
\end{array}
\right.
$$
for all irreducible components 
$C_j\subset C$. Equivalently,  $C\cdot C_j=1$ for every irreducible
component $C_j$.

\begin{proposition}
\label{minus one projective}
Enriques surfaces of type III in minus-one-form are projective.
\end{proposition}

\proof
Let 
$X$ be such a surface.
An invertible 
$\O_S$-module 
$\shL$ descends to 
$X$ if and only if its restriction 
$\shL_C$ lies in   
$\Pic(D)\subset \Pic(C)$. Of course, a necessary condition is that the
numerical class of 
$\shL_C$ lies in 
$\NS(D)\subset\NS(C)$. I claim that this is also sufficient: 
The exact sequence 
$$
0 \lra H^1(X,\O_{X})  \lra H^1(D,\O_D)  \lra H^1(C,\O_C)  
\lra H^2(X,\O_{X})\lra 0
$$
shows that the Jacobians 
$\Pic^0_D$ and $\Pic^0_C$ have the same dimensions. Consequently, the map 
$\Pic_D^0\ra \Pic_C^0$ must be surjective, because its kernel 
$\mu_2$ is infinitesimal. The claim follows.

The divisor $C$ is ample on itself, because $X$ is in minus-one-form.
By the Fujita--Zariski Theorem (\cite{Fujita 1983}, Thm.\ 1.10), the divisor
$C$ is semiample. The corresponding birational contraction 
$S\ra S'$ has an exceptional curve $R\subset S$
disjoint from 
$C$. Choose an $S'$-ample divisor 
$A\in\Div(S)$ supported by $R$. Then 
$tC+A$ is ample for 
$t\gg 0$.
By construction, the class of
$tC+A$  in $\NS(C)$ lies in the subgroup $\NS(D)$. 
Consequently, the ample invertible $\O_S$-module $\O_S(tC+A)$ 
descend to an ample invertible 
$\O_{X}$-module.
\qed

\medskip
In light of this, we seek to put any type III surfaces into minus-one-form. 
To do so, I have to recall
Kulikov's concept of type I and type II modifications \cite{Kulikov 1977}.
Let 
$E\subset D$ be smooth double curve with selfintersection 
$-1$ on each adjacent component.
Blowing-up 
$E$ creates a 
$4$-cycle of 
$(-1)$-curves. Blowing-down the strict transform of 
$E$ yields another Enriques surface of type III called the 
\emph{modification of type II}. The picture is:
\begin{center}
\unitlength1em
\begin{picture}(0,4)
\put(-2,2){\line(-1,1){1}}
\put(-2,2){\line(-1,-1){1}}
\put(-2,2){\line(1,0){4}}
\put(2,2){\line(1,-1){1}}
\put(2,2){\line(1,1){1}}
\put(-.4,1){$\scriptstyle{-1}$}
\put(-.4,2.7){$\scriptstyle{-1}$}
\end{picture}
\hspace{8em}
\begin{picture}(0,4)
\put(-2,2){\vector(1,0){4}}
\put(2,2){\vector(-1,0){4}}
\put(-1.3,2.3){\footnotesize Type II}
\end{picture}
\hspace{6em}
\begin{picture}(0,4)
\put(0,3){\line(-1,1){1}}
\put(0,3){\line(1,1){1}}
\put(0,1){\line(0,1){2}}
\put(0,1){\line(-1,-1){1}}
\put(0,1){\line(1,-1){1}}
\put(-1.6,1.7){$\scriptstyle{-1}$}
\put(.4,1.7){$\scriptstyle{-1}$}
\end{picture}
\end{center}
Now let 
$E\subset S$ be an exceptional curve of the first kind outside 
$C$ corresponding to the dotted lines below. Blowing-up 
$E$ and contracting its   strict transform  yields another Enriques surface 
of type III, called the
\emph{modification of type I}. Here the picture is:
\vspace{-.6em} 
\begin{center}
\unitlength1em
\begin{picture}(0,4)
\put(-2,2){\line(-1,1){1}}
\put(-2,2){\line(-1,-1){1}}
\put(-2,2){\line(1,0){4}}
\put(2,2){\line(1,-1){1}}
\put(2,2){\line(1,1){1}}
\multiput(0,2)(0,.6){3}{\line(0,1){.4}}
\end{picture}
\hspace{6em}
\begin{picture}(0,3)
\put(-2,2){\vector(1,0){4}}
\put(2,2){\vector(-1,0){4}}
\put(-1.3,2.3){\footnotesize Type I}
\end{picture}
\hspace{6em}
\begin{picture}(0,4)
\put(-2,2){\line(-1,1){1}}
\put(-2,2){\line(-1,-1){1}}
\put(-2,2){\line(1,0){4}}
\put(2,2){\line(1,-1){1}}
\put(2,2){\line(1,1){1}}
\multiput(0,2)(0,-.6){3}{\line(0,-1){.4}}
\end{picture}
\end{center}
\vspace{-1em}
 Recall that a normal crossing surface
$X$ is called \emph{d-semistable} if the the sheaf of first-order deformations
$$
\shT^1_X =\shExt^1(\Omega_{X/k}^1,\O_{X})
$$ 
is isomorphic to $\O_D$. Friedman \cite{Friedman 1983}  showed that 
d-semistability is
necessary for the existence of smoothings with smooth total space.

\begin{proposition}
\label{put in form}
For each d-semistable Enriques surface of type III, there is a sequence of 
type I and type II
modifications reaching a d-semistable Enriques surface of type III in 
minus-one-form. 
\end{proposition}

\proof
Miranda and Morrison \cite{Miranda; Morrison 1983} proved this difficult 
result for Kulikov
degenerations of K3 surfaces. Their arguments apply unchanged to our situation.
\qed

\medskip
For later use, we record the following observation.

\begin{lemma}
\label{linear independent}
Suppose 
$X$ is in minus-one-form. Then there are at least two   components 
$S_1,S_2$ so that the irreducible components of 
$C$ contained in $S_1\cup S_2$ have linear independent Hodge classes in the 
$k$-vector space
$H^{1,1}(S)=H^1(S,\Omega_{S/k}^1)$.
\end{lemma}

\proof
First, suppose that the irreducible components of $C\subset S$ are smooth.
Let 
$m>0$ be the number of irreducible components in 
$C$. Then 
$D$ has 
$m/2$ irreducible components and 
$m/3$ singularities. Moreover, the normalization map 
$\tilde{D}\ra D$ has 
$m$ ramification points. Consequently, 
$$
1 = \Pi(X) = m/2+m/3 - m + n = n-m/6.
$$
For each integer 
$r>0$, let 
$n_r$ be the number of components 
$S_i$ so that 
$C\cap S_i$   is an 
$r$-cycle of rational curves. Then 
$n=\sum n_r$ and 
$m=\sum rn_r$, so 
$6 = \sum n_r(6-r)$.
We infer that at least two components, say 
$S_1,S_2$, have a ramification curve with $r\leq 5$ irreducible components.
According to \cite{Miranda; Morrison 1983}, Lemma 11.5, each 
irreducible component 
of $ C$ in $S_1\cup S_2$ meets an exceptional curve of the first kind 
not contained in 
$C$. The assertion follows.

Second, suppose that there is a nodal components $C'\subset C$, 
such that $(C')^2=1$.
If there is another nodal component, we are done immediately,
so let us assume that the other components are smooth.
Then $D$ has $m/2$ irreducible components and 
$(m-1)/3$ singularities. Moreover, the normalization map 
$\tilde{D}\ra D$ has 
$m-1$ ramification points. Now a calculation as above finishes
the proof.
\qed

\medskip
The K3-like covering of a type III Enriques surfaces is a 
K3 surface of type III.
Let me remark that there is no analog of Proposition 
\ref{minus one projective} for such  surfaces:

\begin{proposition}
\label{no embedding}
K3 surfaces of type III with generic gluing map 
do not embed into  smooth separated schemes. 
\end{proposition}

\proof
Let 
$X$ be an  arbitrary K3 surface type III. 
We have to show that   $X$ does not embed into smooth schemes if 
the gluing map $\varphi:C\ra D$ is generic. Here, as usual,  $D$ 
is the singular locus, and 
$C\subset S$ the ramification locus.
Decompose 
$D=D_1\cup\ldots\cup D_m$ into irreducible components. Let 
$k^\times$ act on each
$D_i$ fixing the triple points. This gives 
$\Aut^0(D)=\prod_{i=1}^m k^\times$. We can view 
$T=\Aut^0(D)$ as a parameter space for   gluing morphisms 
$\varphi:C\ra D$ in the following way. For each 
$D_i$, choose an ordering 
$\nu^{-1}(D_i)=C_{i }' \cup C_{i }''  $. Set 
$C'=\cup C_i'$ and $C''=\cup C_i''$. Then each 
$\lambda=(\lambda_1,\ldots,\lambda_m)$ in $T=\Aut^0(D)$ is also an 
automorphism 
$\tilde{\lambda}: C\ra C$, defined by 
$\tilde{\lambda}|_{ C'} =\lambda$ and 
$\tilde{\lambda}|{ C''} =\id$. So each 
$\lambda\in T$ defines a new gluing morphism 
$\varphi_\lambda= \varphi\circ \tilde{\lambda}$. We obtain a 
locally trivial flat family 
$X_\lambda,\lambda\in T$ of type III surfaces.

How does  the Picard group 
$\Pic(X_\lambda)$ jump in this family? 
An invertible 
$\O_S$-module 
$\shL$ descends to 
$X$ if and only if its restriction 
$\shL_C$ lies in   
$\Pic(D)\subset \Pic(C)$. Of course, a necessary condition is that the 
numerical class of 
$\shL_C$ lies in 
$\NS(D)\subset\NS(C)$, in other words, 
$(\shL\cdot C_i')=(\shL\cdot C_i'')$. Suppose this is the case and set 
$d_i=(\shL\cdot C_i')$.
Let 
$\Pic^d(C)\subset \Pic(C)$ be the connected component with
$\shL_C\in\Pic^d(C)$ and 
$\Pic^d(D)\subset\Pic(D)$ the corresponding connected component.
The exact sequence 
$$
0 \lra H^1(D,\O_D)  \lra H^1(C,\O_C)  \lra H^2(X,\O_X)  \lra 0
$$
and $h^2(X,\O_{X})=1$ implies that  there is an exact sequence 
$$
0 \lra \Pic^0(D)  \lra \Pic^0(C)   \lra k^\times  \lra 1.
$$
Hence the inclusion 
$\Pic^d(D)\subset\Pic^d(C)$ has codimension 1. From this we infer that 
the set of all parameters 
$T_\shL\subset T$ for which 
$\shL$ descends to 
$X$ is  either empty, or of codimension 1, or equals 
$T$.

Suppose that 
$\shL_C$ has at least one  degree 
$d_j\neq 0$. Consider parameters
$\lambda=(\lambda_1,\ldots,\lambda_m)$  in $T=\prod_{i=1}^m k^\times$ with 
$\lambda_i=1$ for $i\neq j$. You easily check that  
$\tilde{\lambda}^*(\shL_C)\simeq \shL_C$ holds if and only if 
$\lambda_j\in k^\times$ is a 
$d_j$-th root of unity. It follows that 
$T_\shL\subset T$ is of codimension 1.
Thus the set of all 
$\lambda\in T$ for which some 
$\shL\in\Pic(S)$ with 
$\shL_C\not\equiv \O_C$ descends to 
$X$ is a countable union of codimension 1 subsets.
Consequently, no such 
$\shL$ exists for 
$\lambda\in T$ generic.

Now we argue as follows: Seeking a contradiction, 
suppose that, for generic gluing map $\varphi:C\ra D$, 
there is a closed embedding $X\subset Y$
into a smooth separated scheme
$Y$.  Choose an affine open neighborhood $U\subset Y$ intersecting $D$. 
Then $Y\setminus U$ defines
a Cartier divisor on $Y$, whose restriction to $X$  is 
numerically nontrivial on $D$, contradiction.
\qed

\begin{remark}
\label{divisorial}
Borelli \cite{Borelli 1963} calls a scheme $Y$ is \emph{divisorial} 
if the open subsets of
the form
$Y_s\subset Y$ generate the topology, where 
$s\in\Gamma(Y,\shL)$ ranges over    sections  of   line bundles 
$\shL\in\Pic(Y)$.
The preceding proof shows that it is impossible to embed   
generic K3 surfaces of type III into   divisorial schemes.  
\end{remark}

\section{Fibrations and d-semistability}

Smooth Enriques surface  carry  elliptic or
quasielliptic fibrations (\cite{Cossec; Dolgachev 1989}, Thm.\ 5.7.1).
Do semistable  Enriques surface admit a genus one fibration? 
In this section we shall find partial answers for   type II surfaces.
First, we observe that there is always  an ample line bundle:

\begin{proposition}
\label{type II projective}
 Enriques surface of type II are projective.
\end{proposition}

\proof
Let $X$ be such a surface.
To construct an ample invertible 
$\O_{X}$-module, it suffices to find an ample invertible 
$\O_S$-module $\shL$ whose restriction 
$\shL_C$ lies in the image of 
$\Pic(D)$.
Let us do the case that 
$X$ is nonsimple and reducible (the other cases are similar).
We use the notation from Theorem \ref{nonsimple type II}.
 Start with the rational component 
$S_1$ and choose an ample invertible 
$\O_{S_1}$-module 
$\shL_1$. Next, choose an ample invertible 
$\O_{S_2}$-module 
$\shL_2$. Passing to suitable multiples if necessary, we can assume that 
$(\shL_1\cdot C_1)=(\shL_2\cdot C_2'')$, where 
$C_2''\subset C_2$ is the connected component with 
$\varphi(C_1)=\varphi(C_2'')$ . Since 
$\Pic^0(S_2)\ra\Pic^0(C_2'')$ is surjective, we can modify 
$\shL_2$ by a numerically trivial sheaf so that 
$\shL_1|_{ C_1}\simeq \shL_2|_{ C_2''}$ with respect to the gluing morphisms. 

Now proceed by induction: This gives a sequence 
$\shL_1,\ldots,\shL_n$ of ample invertible 
$\O_{S_i}$-modules which coincides on the overlaps.
Let 
$D_n\subset D$ be the image of the last ramification curve
$C_n'\subset C_n$. Then $C_n'\ra D_n$ is a $\ZZ/2\ZZ$ double covering, 
such that
$\Pic(D_n)\ra\Pic(C_n')$ has a  cokernel $\ZZ/2\ZZ$.
 Replacing the
$\shL=\shL_1\cup\ldots \shL_n$ by a suitable multiple if necessary, we can assume that 
$\shL$ is an  invertible
$\O_S$-module whose restriction 
$\shL_C$ lies in the image of 
$\Pic(D)$.
\qed

\medskip
In the following, 
$X$ will be a  Enriques surface of type II. Then the sheaf 
$\shT_X^1 =\shExt^1(\Omega_{X/k},\O_{X})$ 
 is  dual to the norm 
$N_{C/D}(\shN)$, where 
$\shN=\shN_{C/S}$ is the conormal bundle. Let me state the  following fact:

\begin{lemma}
\label{d-semistable double}
Suppose $X$ is an  Enriques surface of  type II. Let $D'\subset D$ be a 
connected component
with connected preimage
$C'\subset C$. Then 
$\shT^1_X|_{D'}\simeq \O_{D'}$ holds if and only if the conormal bundle 
$\shN_{C'/S}\in\Pic(C')$ is 2-torsion.
\end{lemma}

\proof
Set $\shN=\shN_{C'/S}$. We  have an exact sequence 
$$
0 \lra D(\ZZ/2\ZZ)  \lra \Pic^0_{D'}  \lra \Pic^0_{C'}  \lra 0,
$$
so $\Pic^0(D')\ra \Pic^0(C')$ is bijective. Suppose 
$\shN^{\otimes 2}\simeq \O_{C'}$. Choose an invertible 
$\O_{D'}$-module 
$\shL$ with 
$\shN=\nu^*(\shL)$. Then 
$\N_{C'/D'}(\shN)=\shL^{\otimes 2}$ is trivial, so 
$\shT^1_X|_{D'}\simeq \O_{D'}$. The converse is similar.
\qed

\medskip
The condition of d-semistability does not restrict the number of 
irreducible components of 
$X$. Rather, it says something about the minimal model of the normalization:

\begin{proposition}
\label{minimal model}
Suppose 
$X$ is   d-semistable and of type II. Let 
$\bar{S}$ be a minimal model of 
$S$. Then the corresponding contraction
$h:S\ra \bar{S}$ is a sequence of 8 or 9 blowing-ups.
\end{proposition}

\proof
The number $m\geq 0$ of blowing-ups in 
$h:S\ra \bar{S}$ is given by  $K_S^2=K_{\bar{S}}^2-m$.
On the other hand,    d-semistability implies
$ K_S^2=C^2=0$. 
According to Theorem \ref{nonsimple type II} and Theorem\ref{simple type II}, 
there is precisely one rational components
$\bar{S}_1\subset \bar{S}$, which has 
$K_{\bar{S}_1}^2=8$ or 
$K_{\bar{S}_1}^2=9$. Moreover, the other minimal models have 
$K_{\bar{S_i}}^2=0$. This proves the assertion.
\qed
 
\medskip 
Regular but nonsmooth curves of arithmetic genus one are called \emph{quasielliptic}.
They do not exist  over algebraically closed fields, but may occur   over   function fields:

\begin{proposition}
\label{type II genus 1 fibrations}
Suppose that
$X$ is  a nonsimple d-semistable    Enriques surface of type II.
Then there is a proper morphism 
$g:X\ra \PP^1$ with 
$\O_{\PP^1}\ra g_*(\O_{X})$   bijective, so 
that the generic fiber $X_\eta$ is an elliptic or
quasielliptic curve.
\end{proposition}

\proof
First, assume    that 
$X$ is irreducible. According to Theorem \ref{nonsimple type II}, 
the normalization 
$S$ is rational, and the ramification locus
$C\subset S$ is an elliptic curve. By Lemma \ref{d-semistable double}, the conormal sheaf 
$\shN=\O_C(-C)$ is 2-torsion.
Assume that 
$\shN$ is trivial. Then the exact sequence 
$$
0 \lra \O_{X}  \lra \O_{X}(C)  \lra \O_C   \lra 0
$$
yields an exact sequence 
$$
0 \lra H^0(S,\O_{S})  \lra H^0(S,\O_{S}(C))   \lra H^0(C,\O_{C})   \lra 0,
$$
hence 
$C$ is base-point-free and defines the desired fibration.
Now assume that 
$\shN$ has order 2. The exact sequence 
$$
0 \lra \O_{S}  \lra \O_{S}(2C)  \lra \O_{2C}(2C)   \lra 0
$$
yields an exact sequence 
$$
0 \lra H^0(S,\O_{S})  \lra H^0(S,\O_{S}(2C))   
\lra H^0(2C,\O_{2C}(2C))   \lra 0.
$$
Moreover, the exact sequence
$$
0 \lra \\O_C(C) \lra \O_{2C}(2C)  \lra \O_{ C}  \lra 0
$$
implies that 
$H^0( \O_{2C}(2C))\ra H^0(\O_C)$ is bijective, hence 
$2C$ is base-point-free.
In both cases we obtain a genus one fibration 
$S\ra \PP^1$, which clearly descends to a genus one fibration 
$g:X\ra \PP^1$.

Next, assume that     
$X$ is not  irreducible. Using the notation from Theorem \ref{nonsimple type II}, 
we shall consider the last component  $S_n\subset S$ and its ruling $f_n:S_n\ra B_n$.
Let $C'\subset C_n$ be the connected component double covering its image $D'\subset D$.
By d-semistability,   $\shN=\shN_{C'/S}$ is 2-torsion.
Setting $\shE=(f_n)_*(\O_{S_n}(C'))$ and $\shL=(f_n)_*(\shN)$, we obtain an exact sequence 
$$
0 \lra  \O_{B_n}  \lra \shE   \lra  \shL \lra 0.
$$
Note that 
$\bar{S}_n=\PP(\shE)$ gives a minimal model 
$h:S_n\ra\bar{S}_n$.
Suppose that 
$\shN$ is trivial. Since the other ramification curve $C''\subset C_n$ has image
$h(C_n'')$   disjoint from   
$h(C')$, the preceding extension splits, and 
$\bar{S}_n=B_n\times\PP^1$. So we obtain an elliptic structure 
$S_n\ra \PP^1$. 

Now suppose that 
$\shN$ has order 2. The corresponding inclusion 
$\ZZ/2\ZZ\subset \Pic_{B_n} $ defines a principle homogeneous 
$\mu_2$-space over
$B_n$ on which 
$\shN$ becomes trivial. It follows that 
$\Fr^*(\shN)$ is trivial. We infer that there is an elliptic  structure on 
$\PP(\Fr^*(\shE))$ which descends to a genus one fibration on 
$\bar{S}_n=\PP(\shE)$. In both cases, the fibrations induce a morphism 
$X\ra \PP^1$. Note that the initial components 
$S_1,\ldots, S_{n-1}$ are mapped to points.
\qed

\medskip
Here is another kind of fibration. Suppose that 
$X$ is of type II as described in Theorem \ref{nonsimple type II} or \ref{simple type II}.
Additionally, assume  that the rational component 
$S_1\subset S$ is not 
$\PP^2$. Choose a ruling 
$f_1:S_1\ra \PP^1$. Together with the elliptic rulings 
$f_i:S_i\ra B_i$ for 
$2\leq i\leq n$, this defines a ruling 
$f:S\ra B$ over 
$B=\PP^1\cup B_2\cup\ldots\cup B_n$. We seek to descend this fibration on $S$ to  a fibration on $X$ mapping to a
finite quotient of the curve $B$.

\begin{proposition}
\label{second fibration on type II}
With the preceding assumptions, the fibration 
$f:S\ra B$ descends to a fibration 
$g:X\ra\PP^1$. The generic fiber is a curve (possibly nonreduced)
of arithmetic genus 
$p_a(X_\eta)=1$.
\end{proposition}

\proof
Let 
$E$ be the elliptic curve isomorphic to the components of 
$C$. The induced projection 
$f_1:C_1\ra \PP^1$ defines an involution 
$\iota:E\ra E$, hence a subgroup 
$\ZZ/2\ZZ\subset\Aut(E)$.
There is another subgroup scheme 
$H\subset \Aut(E)$ of length two:
If 
$X$ is simple, we define
$H$  by the double covering 
$C_n\ra B_n$; if 
$X$ is nonsimple, 
$H=\ZZ/2\ZZ$ is defined by the part of the gluing morphism that is 
responsible for the nonsmooth component 
$X_n$.
Let 
$G\subset \Aut(E)$ be the subgroup scheme of length 4 generated by 
$H$ and the involution 
$\iota:E\ra E$. Then 
$E/G\simeq \PP^1$.  It is easy to see that 
$f:S\ra B$ induces a fibration 
$g:X\ra E/G$.

It remains to determine the generic fiber.
First, suppose that 
$X$ is nonsimple.  Then 
$X_\eta$ is a cycle of rational curves with 
$2+4(n-1) = 4n-2$ irreducible components. Here each elliptic ruled component 
$S_i$ contains 4 irreducible components of 
$X_\eta$.
Second, suppose that 
$X$ is simple. Then
$X_\eta$ is nonreduced: it is a string of the form 
$X_\eta = \PP^1 + 2\PP^1 +\ldots+2\PP^1 + \PP^1$ with 
$1+2(n-1) = 2n-1$ irreducible components. Here the reduced part 
$\PP^1+\PP^1$ lies on the last irreducible component 
$X_n\subset X$.
From this, you  immediately get 
$p_a(X_\eta)=1$.
\qed

\section{Using  log structures to construct formal smoothings}
\label{Using  log structures to construct formal smoothings}

In this section, the task is to find \emph{formal} deformations of 
d-semistable Enriques surfaces.
This will be a major step towards the construction of algebraic deformations, 
which eventually leads to 
smoothings and liftings. 
A direct approach would be to study a versal deformation 
$\foX\ra \Spf(R)$. The problem, however, is that $X=\foX_0$ has many 
locally trivial deformations, 
which
are irrelevant for our purposes, and this fact obscures the structure of 
$\Spec(R)$. 

A natural  way to avoid such problems is to use \emph{logarithmic structures}.
The underlying idea is to enlarge the category of schemes so that d-semistable
schemes can be considered as smooth schemes. Let my recall the 
fundamental definition.

\begin{definition} 
(Kato)
A \emph{log scheme} 
$Y^\dagger = (Y,\shM_Y,\alpha)$ comprises a scheme 
$Y$, a sheaf of monoids
$\shM_Y$,  and  a homomorphism  
$\alpha:\shM_Y\ra \O_{Y}$ so that  the map
$\alpha^{-1}(\O_{X}^\times)\ra \O_{X}^\times$  is bijective.
\end{definition}

Here 
$\shM_Y$ is a sheaf   in the   \'etale topology, and 
$\alpha:\shM_Y\ra \O_{Y}$ is a   homomorphism with respect to 
multiplication in 
$\O_{X}$. You  find more about log structures in K.~Kato's 
fundamental paper \cite{Kato 1989} and
Illusie's survey article
\cite{Illusie 1994}.
Suppose   
$X$ is  a d-semistable   normal crossing   surface.
We shall use log structures on $X$, on its normalization 
$S$, and on the ground field $k$.

Let us start with the ground field.
Set 
$\shM_k=\NN\oplus k^\times$. The function 
$\alpha:\shM\ra k$ with 
$\alpha(0,\lambda)=\lambda$ and 
$\alpha(n,\lambda)=0$ for 
$n\neq 0$ defines a log structure 
$k^\dagger$ called the \emph{standard log structure}. The corresponding  
log scheme 
$\Spec(k^\dagger)$ is  nicknamed the `punctured point'.

Next, consider the normalization 
$S$ of $X$. Let 
$j:U\ra S$ be the complement of the ramification locus 
$C\subset S$ for normalization. Set 
$\shM_S=\O_S \cap j^*(\O_{U}^\times)$. The canonical map 
$\alpha:\shM_S\ra\O_S$ yields a log structure 
$S^\dagger$.

Finally, we come to the normal crossing surface $X$.
Locally, there is a closed embedding $i:X\ra \AA^3$ as  a union of 
coordinate hyperplanes.  Pulling
back the log structure
$\O_{\AA^3} \cap \O_{\AA^3\setminus X}^\times $  defines local  
log structures on $X$. Of
course, these local log structures might not be compatible. However,
F.~Kato \cite{Kato 1996}, Theorem
11.7 showed that a compatible choice is possible if 
$X$ is d-semistable. The corresponding global log structures are 
called of \emph{semistable type}.
Fix such  a log structure of semistable type $X^\dagger$.
By \cite{Kato 1996}, Example 4.7,   
the structure morphism  induces a log smooth morphism
$X^\dagger\ra  \Spec(k^\dagger)$.
This means that   the lifting criterion for smoothness holds
in the category of log schemes.

We seek to extend log structures of semistable type $X^\dagger$ over 
local Artin rings.
Let $W$ be a complete discrete valuation ring of mixed characteristic with 
residue field $k$.
(Choose $W\subset W(k)$ as the Cohen subring of the ring of Witt vectors, 
as in Bourbaki 
\cite{AC 8-9}, Chap.~IX.)

To keep track of the log structures, we introduce  a formal variable $T$. The 
power series ring
$W[[T]]^\dagger$ is endowed with the log structure given by $\NN\ra W[[T]]$, 
$n\mapsto T^n$. Let 
$(\Art_{W[[T]]})$ be the category of Artin local 
$W[[T]]$-algebras   with residue field  
$k$. Each such Artin ring
$A$ inherits a   log structure
$A^\dagger$    from the log structure 
$W[[T]]^\dagger$. For example, the residue field 
$k$ is endowed with the standard log structure.

Let 
$\LD_{X^\dagger}(A)$ be the set, modulo isomorphism, of   log smooth morphisms 
$Y^\dagger\ra \Spec(A^\dagger)$  extending the log smooth morphism
$X^\dagger\ra\Spec(k^\dagger)$. We call the corresponding  functor 
$\LD_{X^\dagger}:(\Art_{W[[T]]})\ra (\Set)$ the 
\emph{log deformation functor} of 
$X^\dagger$. The following is   the main technical results of this paper:

\begin{theorem}
\label{functor smooth}
Suppose  $X$ is  a nonsmooth d-semistable   Enriques surface.
If $X$ is classical or ordinary, then the log deformation functor
$\LD_{X^\dagger}:(\Art_{W[[T]]})\ra (\Set)$  is formally smooth.
\end{theorem}

\proof
Formal smoothness of 
$\LD_{X^\dagger}$ means: For each log smooth deformation 
$Y^\dagger/A^\dagger$ and each extension 
$B\ra A$ with square-zero ideal 
$I\subset B$, it is possible to extend 
$Y^\dagger$ over 
$B^\dagger$.
According to F.~Kato \cite{Kato 1996},
Proposition 8.6, the obstruction lies in 
$H^2(X,\shHom(\Omega^1_{X^\dagger/k^\dagger},\O_{X}))\otimes_A I$.  Here 
$\Omega^1_{X^\dagger/k^\dagger}$ is the sheaf of 
\emph{log differentials} (see \cite{Kato 1996}, sect.~5), 
which is locally free.
It suffices to check 
$H^0(X,\Omega^1_{X^\dagger/k^\dagger}\otimes\omega_X)=0$.

Do proceed we have to relate various differentials. Luckily, this was done in
Friedman's paper \cite{Friedman 1983},  where the sheaf 
$\Omega^1_{X^\dagger/k^\dagger}$ appears under the name 
$\Lambda_X^1$. He constructed a commutative diagram of 
$\O_{X}$-modules
\begin{equation}
\label{diagram}
\begin{CD}
@. 0 @. 0 @. 0\\
@. @VVV @VVV @VVV\\
0     @>>>  (\Omega_{X/k}^1)^{\vee\vee} @>>> \Omega_{X^\dagger/k^\dagger}^1 
@>>> 
\O_{\tilde{D}} @>>> \O_T  @>>> 0  \\
@. @VVV     @VVV @VVV\\
0 @>>> \Omega_{S/k}^1 @>>> \Omega_{S^\dagger/k}^1 @>>> \O_{\tilde{C}}    
@>>>  0\\
@. @VVV\\\
@. \Omega_{\tilde{D}/k}^1\\
@. @VVV\\
@. 0
\end{CD}
\end{equation}
with exact rows and columns.
Here 
$\O_{\tilde{D}}\ra\O_D$ and 
$\O_{\tilde{C}}\ra\O_C$ are the normalizations, and  
$T\subset X$ is the set of triple points.
The second row displays $\Omega^1_{S^\dagger/k}$ as an inverse
elementary transformation of $\Omega_{S/k}^1$.
The map 
$\Omega_{S/k}^1\ra \Omega_{\tilde{D}/k}^1$ is a \v{C}ech boundary operator 
explained in \cite{Friedman 1983}, p.~77.  We shall treat three cases.

\emph{(i) The case that $X$ is ordinary and of  type II.}
There are no triple points, so 
$D=\tilde{D}$ and 
$C=\tilde{C}$.  The upper row of Diagram (\ref{diagram}) 
gives an exact sequence 
$$
H^0(X,(\Omega_{X/k}^1)^{\vee\vee})  \lra H^0(X,\Omega^1_{X^\dagger/k^\dagger})
\lra H^0(D,\O_D)\lra H^1(X,(\Omega_{X/k}^1)^{\vee\vee}).
$$
First, I claim that 
$H^0(X,(\Omega_{X/k}^1)^{\vee\vee})=0$. Indeed: The left column of 
Diagram (\ref{diagram}) gives an
exact sequence 
$$
0 \lra H^0(X,(\Omega_{X/k}^1)^{\vee\vee})   \lra H^0(S,\Omega_S^1)  \lra 
H^0(D, \Omega_D^1).
$$
Using the notation from Theorem \ref{simple type II}, we have 
$H^{1,0}(S_1)=0$ because the component
$S_1$ is rational. Moreover, for 
$2\leq i\leq n$, the maps in 
$$
H^{1,0}(B_i)  \lra H^{1,0}(S_i)   \lra   H^{1,0}(C\cap S_i) 
$$
are injective because the projections 
$C\cap S_i\ra B_i$ are \'etale. You easily infer that 
$H^{1,0}(S) \ra H^{1,0}(D)$ is injective, so 
$H^0(X,(\Omega_{X/k}^1)^{\vee\vee})=0$.

Second, I claim that the boundary map
$H^0(\O_D) \ra H^1(X,(\Omega_{X/k}^1)^{\vee\vee})$ is injective. 
To see this, consider the commutative
diagram 
$$
\begin{CD}
H^0(D,\O_D) @>>> H^1( X,(\Omega_{X/k}^1)^{\vee\vee})  \\
@VVV     @VVV\\
H^0(C,\O_C) @>>> H^1(S,\Omega_{S/k}^1).  
\end{CD}
$$
It suffices to check that the composition 
$H^0(D,\O_D)\ra H^{1,1}(S)$ is injective. A direct manipulation with 
cocycles shows that the boundary
map 
$H^0(C,\O_{C})\ra H^{1,1}(S)$ maps a section of 
$\O_{C}$ to the  Hodge class of its support. Recall that the Hodge-class-map
$\Pic(S)\ra H^{1,1}(S)$ is induced by logarithmic derivation 
$$
\dlog:\O_S^\times\lra \Omega_X^1,\quad s\mapsto \dlog(s)=ds/s.
$$
Let 
$C_j\subset C$ be the irreducible components, and assume that some divisor 
$\sum\lambda_jC_j$ coming from 
$H^0(D,\shF)$ has zero Hodge class.
By Proposition \ref{minimal model}, there is an exceptional curve of the 
first kind 
$E\subset S$, say with 
$E\cdot C_j=1$. Then 
$\lambda_j=0$. Let 
$C_k$ be another component. If 
$\varphi(C_k)=\varphi(C_j)$ is a double curve, then 
$\lambda_k=\lambda_j$ is zero as well. If 
$C_j,C_k$ lie in the same ruled component 
$S_i$, the ruling 
$S_i\ra B_i$ implies that 
$\lambda_k=\lambda_j=0$.
Inductively, we conclude that  
all multiplicities in 
$\sum\lambda_jC_j$ vanish.
It follows that 
$H^0(X,\Omega^1_{X^\dagger/k^\dagger})=0$, and we conclude that the 
log deformation functor is
formally smooth.

\emph{(ii) The case that X is ordinary of type III.}
According to Theorem \ref{type III}, the components 
$S_i$ are rational surface, so both 
$H^0(X,(\Omega^1_X)^{\vee\vee})\subset H^{1,0}(S)$ must vanish.  
Consider the commutative diagram 
$$
\begin{CD}
H^0(D,\shF)     @>>>  H^1(X,(\Omega^1_{X/k})^{\vee\vee})  
\\ @VVV @VVV  \\
H^0(\tilde{C},\O_{\tilde{C}})     @>>>   H^1(S,\Omega^1_{S/k}).
\end{CD}
$$
with $\shF=\Omega^1_{X^\dagger/k^\dagger} / ( \Omega^1_{X/k})^{\vee\vee}$. 
We have to check that
$H^0(D,\shF)\ra  H^1(X,(\Omega^1_{X/k})^{\vee\vee})$ is injective.
 It suffices to show that the composition
$H^0(D,\shF)\ra H^{1,1}(S)$ is injective. 
Let 
$C_j\subset C$ be the irreducible components and consider a nonzero divisor
$\sum\lambda_j C_j$ coming from 
$H^0(D,\shF)$. This means that 
$\lambda_j=\lambda_k$ if the intersection
$\varphi(C_j)\cap\varphi(C_k)$ is a double curve, and 
$\lambda_j+\lambda_k=\lambda_l$ if the intersection
$\varphi(C_j)\cap\varphi(C_k)\cap\varphi(C_l)$ is a triple point.

 Seeking a contradiction, we assume that 
$\sum\lambda_j C_j$   has zero Hodge class.
By Proposition \ref{put in form}, suitable modifications of type I and 
type II put 
$X$ into minus-one-form 
$X'$. Let 
$C'\subset S'$ be the corresponding ramification curve, and 
$C'_j\subset C'$ be its irreducible components. The modifications give 
a canonical bijection
between the 
$C_j$ and the 
$C_j'$, and the nonzero divisor 
$\sum\lambda_jC'_j$   has zero Hodge class in $H^{1,1}(S')$. 
According to Proposition 
\ref{linear independent}, there is a component 
$S_0'\subset S'$ with 
$\lambda_j=0$ for all
$C'_j\subset S'_0$. By the triple point condition, the multiplicities near 
$S'_0$ look like:
\begin{center}
\unitlength1em
\begin{picture}(0,4)
\put(-2,2){\line(-1,1){2}}
\put(-2,2){\line(-1,-1){2}}
\put(-2,2){\line(1,0){4}}
\put(2,2){\line(1,-1){2}}
\put(2,2){\line(1,1){2}}
\put(-6,1.8){$S'_0$}
\put(-2.5,3){$0$}
\put(-2.5,0){$0$}
\put(-.2,1){$0$}
\put(2,3){$$}
\put(2,0){$$}
\end{picture}
\end{center}
Using that 
$\sum\lambda_jC_j'$ has zero Hodge class and that
$X'$ is in minus-one-form, we inductively infer  that 
$\lambda_j=0$ for all 
$C_j$, contradiction.  Again the log deformation functor is
formally smooth.  

\emph{(iii) The case that X is classical.}
For the sake of simplicity, I only do the case that 
$X=X_1\cup X_2$ has 
$n=2$ irreducible components. By Theorem \ref{simple type II}, the double curve
$D=X_1\cap X_2$ is irreducible and 
$C$ comprises two irreducible components $C_i\subset S_i$.
By Proposition \ref{ramification on elliptic components}, we have
$\nu^*(K_X)|_{S_1}=0$, whereas 
$\nu^*(K_X)|_{S_2}$ has order two. You easily deduce that both
$H^0(X,(\Omega^1_X)^{\vee\vee}\otimes\omega_X)\subset 
H^0(S,\Omega^1_{S/k}\otimes\omega_X)$ vanish. 

It remains to check that the boundary map
$H^0(\O_D) \ra H^1(X,(\Omega_{X/k}^1)^{\vee\vee}\otimes\omega_X)$ 
is injective. 
Seeking a contradiction, we
suppose   that this map  is zero. Then the map
$H^1(X,(\Omega_{X/k}^1)^{\vee\vee}\otimes\omega_X)\ra
H^1(X,\Omega^1_{X^\dagger/k^\dagger}\otimes\omega_X)$ is injective,  
so the dual map 
$$
H^1(X,\shHom(\Omega^1_{X^\dagger/k\dagger},\O_{X}))\lra 
H^1(X,\shHom(\Omega^1_{X/k},\O_{X}))
$$
is surjective. The group on the left is nothing but the tangential space 
$\LD_{X^\dagger}(k[\epsilon])$ for all locally trivial log deformations. 
The group on the right 
is the tangential space 
$H^1(X,\Theta_X)$ for all locally trivial deformations.
The preceding surjection means that each first-order 
locally trivial deformation
can be endowed with a log structure.

But this is absurd:
By Proposition \ref{minimal model}, there is an exceptional curve of the 
first kind 
$E\subset S$. Let 
$h:S\ra\bar{S}$ be its contraction. Moving  the center 
$\bar{s}=h(E)$ in 
$h(C)\subset \bar{S}$ over the dual numbers 
$k[\epsilon]$ destroys d-semistability
for the corresponding first-order deformation of 
$X$. This gives an element in 
$H^1(X,\Theta_X)$ not in the image of 
$\LD_{X^\dagger}(k[\epsilon])$, contradiction.
Again, we conclude that the log deformation functor is formally smooth.
\qed

\section{Supersingularity and obstructions}

In contrast to the classical and ordinary case, there are obstructions for 
supersingular Enriques
surfaces.

\begin{theorem}
\label{supersingular functor nonsmooth}
The log deformation functor 
$\LD_{X^\dagger}:(\Art_{W[[T]]})\ra (\Set)$ of a nonsmooth d-semistable 
supersingular
Enriques surface  is not formally smooth.
\end{theorem}

This needs some preparations.  First, let me recall 
the classification of group schemes of order two over an 
arbitrary  ground ring 
$A$  of characteristic 
$p=2$. Fix two elements 
$a,b\in A$ with 
$ab=0$. Set 
$$
\Lambda_a=A[T |T^2-aT]\quadand
\mu_b(T)= T\otimes 1+1\otimes T + b T\otimes T.
$$
Then 
$G_{a,b} =\Spec(\Lambda_a)$ is a commutative group scheme of length two, with 
group law defined by 
$\mu_b:\Lambda_a\ra\Lambda_a\otimes\Lambda_a$. 
The corresponding group-valued functor on the category of  $A$-algebras is
$G_{a,b}(R)=\left\{  r\in R\mid r^2=ar \right\}$
with group law 
$$
r_1*r_2 = r_1+r_2 + br_1r_2.
$$ 
For example, 
$G_{0,0}=\alpha_2$, and 
$G_{1,0}=\ZZ/2\ZZ$, and 
$G_{0,1}=\mu_2$. 
There is an isomorphism 
$G_{a,b}\simeq G_{a',b'}$ if and only if there is a unit 
$\gamma\in A^\times$ with 
$a'=\gamma a$ and 
$b'=\gamma^{-1}b$. Swapping the indices gives Cartier duality: 
$D(G_{a,b})=G_{b,a}$.
Note that the augmentation ideal 
$T\Lambda\subset \Lambda$ is a free  
$A$-module of rank 1. 

\begin{lemma}
\label{group schemes}
Let  
$G=\Spec(\Lambda)$ be a commutative affine group $A$-scheme 
whose augmentation ideal 
is a free   
$A$-module of rank 1. Then 
$G$ is isomorphic to 
$G_{a,b}$ for some 
$a,b\in A$ with 
$ab=0$.
\end{lemma}

This is a  straightforward calculation. You   find it  in 
\cite{Tate 1995}, Example 3.2.
We shall apply this result in the following situation. Set $A=k[[t]]$, and let 
$\shE\subset \PP^2_A$ be the relative elliptic curve defined by the 
Weierstrass equation in Deuring normal form
$y^2+t^2xy+y=x^3$. The closed fiber has  
$j(\shE_\sigma)=0$, whereas 
$j(\shE_\eta)= t^{33}/(t^6-1)$. Let 
$\shG$ be the kernel   over 
$A$ of the homomorphism
$2:\shE\ra\shE$.

\begin{lemma}
\label{kernel}
The group $A$-scheme 
$\shG$ is isomorphic to 
$G_{t,0}$.
\end{lemma}

\proof
The closed fiber is 
$\shG_\sigma=G_{0,0}$. Using the group law on the generic fiber 
\cite{Silverman 1986}, p.~58, you see
that 
$\shG_\eta=\ZZ/2\ZZ$, generated by the rational point 
$(t^{-2},t^{-3})\in\shE(\eta)$. 
(This is why we use the coefficient $t^2$  instead of $t$ in the 
Weierstrass equation.
Otherwise $2:\shE\ra\shE$ would be a nontrivial $\shG$-torsors.)
According to Lemma \ref{group schemes}, we have 
$\shG=G_{t^m,0}$ for some integer 
$m\geq 1$. Note that $m$ depends only on the underlying scheme
structure of $\shG$; in fact, it is the order of contact for the
closures of the two points in $\shG(\eta)$.
The point $(t^{-2},t^{-3})\in\shG(\eta)$ has homogeneous coordinates 
$(t:1:t^3)1152$, and we conclude $m=1$.
\qed

\medskip
Next, we use the flat family 
$\shE\ra\Spec(A)$ of elliptic curves to construct flat families of type 
II Enriques surfaces.

\begin{lemma}
\label{supersingular deforms}
Suppose $X$ is a supersingular d-semistable Enriques surface of type II.
Then there are two
projective $A$-deformation 
$\foX_1,\foX_2$ of 
$X$ so that 
$\foX_2\otimes\k(\eta)$ is an ordinary  d-semistable Enriques surface of 
type II, whereas 
$\foX_1\otimes\k(\eta)$ is a classical d-semistable  Enriques surface of 
type II . 
\end{lemma}

\proof
We shall use the notation from Theorem \ref{simple type II}.
Setting 
$\shB_i=\shE$, we obtain   deformations of the elliptic curves defined 
by the rulings 
$f_i:S_i\ra B_i$. Next, we deform the ramification locus 
$C\subset S$. Set 
$\shC_i=\shE$ for 
$i=1$ and $i=n$,  and 
$\shC_i=\shE\coprod \shE$ for 
$2\leq i\leq n-1$. Let 
$f_n:\shC_n\ra\shB_n$ be the double covering defined by 
$2:\shE\ra\shE$. For 
$2\leq i\leq n-1$, let 
$f_i:\shC_i\ra\shB_i$ be the disjoint union of  the identity  on
$\shE$.
Now the family of geometrically  ruled surfaces
$\bar{\shS}_i=\PP((f_i)_*(\O_{\shC_i}))$ is a deformation of the 
minimal models 
$\bar{S}_i$. Lifting the centers for the blowing-up in 
$S_i\ra\bar{S}_i$, you obtain a  deformation 
$\shS$ of   
$S$. Deforming the gluing map 
$\varphi:C\ra D$ and making a gluing over $A$, we obtain a flat family 
$\foX\ra \Spec(A)$ of type II Enriques surfaces with special fiber
$\foX_\sigma = X$, such that 
$\foX_\eta$ is an ordinary Enriques surface of type II.

To obtain classical surfaces as generic fibers, we have to replace 
the projection 
$f_n:\shC_n\ra \shB_n$ defined in the previous construction as  
$2:\shE\ra\shE$ with the dual projection 
$2^*: \Pic^0_{\shE/A}\ra \Pic^0_{\shE/A} $, which has kernel 
$G_{0,t}$.
\qed

\medskip
\noindent \emph{Proof of Theorem \ref{supersingular functor nonsmooth}.}  
Set 
$A_1=k[t_1]/( t_1^2)$, and let 
$X_{1}$ be a first-order deformation of 
$X$ towards ordinary Enriques surfaces obtained from Lemma 
\ref{supersingular deforms} by restriction.
Then the 
relative Picard scheme is
$\Pic^\tau_{X_{1}/A_1}=G_{0,t_1}$.
Similarly, set 
$A_2=k[t_2]/( t_2^2)$, and let 
$X_{2}$ be the first-order deformation of 
$X$ towards classical Enriques surfaces.
Then 
$\Pic^\tau_{X_{2}/A_2}=G_{t_2,0}$.

Seeking a contradiction, we assume that the log deformation functor
$\LD_{X^\dagger}$ is formally smooth. Setting 
$A=k[t_1,t_2]/(t_1^2, t_2^2)$, we can find a    deformation  
$\foX$ over $A$ extending both 
$X_{1}$ and 
$X_{2}$. Consider the relative Picard scheme 
$G=\Pic^\tau_{\foX/A}$. Its closed fiber is 
$G\otimes k= G_{0,0}$. By Nakayama's Lemma,
$G$ is a closed subgroup scheme of  
$G_{a,b}$ for some 
$a,b\in A$ with 
$ab=0$.

By construction,
$G_{a,b}\otimes A_1= G_{0,t_1} $.  Calculating modulo $t_2$, we obtain 
$a\equiv 0$ and 
$b\equiv \lambda_1t_1$ for some unit 
$\lambda_1\in k$.
Similarly, 
$G_{a,b}\otimes A_2= G_{t_2,0} $. Calculating modulo 
$t_1$,  we have  
$b\equiv 0$, and
$a\equiv \lambda_2 t_2$  for some unit 
$\lambda_2\in k$.
This gives
$a=\lambda_2 t_2$ and 
$b=\lambda_1t_1$, contradicting
$ab=0$.
\qed

\section{Smoothings and liftings}

Let 
$X$ be a proper scheme over a field 
$k$ of characteristic 
$p=2$. We shall say that 
$X$ is \emph{smoothable} if there is an integral local noetherian ring 
$A$ with residue field 
$k$, together with a proper flat 
$A$-scheme 
$Y$ with closed fiber
$Y_\sigma=X$ and smooth generic fiber. Furthermore, we say that 
$X$ \emph{lifts to characteristic zero} if there is an integral 
local noetherian ring
$A$ of mixed characteristic with residue field 
$k$, together with a  proper flat 
$A$-scheme 
$Y$ with closed fiber 
$Y_\sigma=X$. Now we come to the main result of this paper.

\begin{theorem}
\label{smoothable and lift}
Nonsmooth d-semistable Enriques surfaces with normal crossings are smoothable
and lift to characteristic zero.
\end{theorem}

\proof
Fix such an Enriques surface 
$X$ and choose a log structure 
$X^\dagger$ of semistable type. Let 
$W$ be a discrete valuation ring  of mixed characteristic with residue field 
$k$. According to \cite{Kato 1996}, Theorem 8.7, 
the corresponding log deformation
functor 
$\LD_{X^\dagger}:(\Art_{W[[T]]})\ra  (\Set)$ admits a hull 
$h:\foX^\dagger\ra\Spf(R^\dagger)$. Here 
$R$ is a complete local  noetherian
$W[[T]]$-algebra, and 
$\foX^\dagger$ is a proper   formal log smooth
$R^\dagger$-scheme.
We shall proceed in five steps.

\textit{Step (i): $X$ is classical.}
According to Proposition \ref{type II projective}, there is an ample 
$\O_{X}$-module 
$\shL$. Let 
$X_m\subset\foX$ be the 
$m$-th order infinitesimal neighborhood of the closed fiber $X\subset \foX$. 
The exact sequence 
$$
\Pic(X_{m+1})   \lra \Pic(X_m)   \lra H^2(X,\O_{X})
$$
together with 
$H^2(X,\O_{X})=0$ ensures that 
$\shL$ extends to an invertible 
$\O_{\foX}$-module.  
Hence Grothendieck's Existence Theorem (\cite{EGA IIIa}, Thm.\ 5.4.5) 
applies, and
we conclude that 
$\foX$ is algebraizable. In other words, there is a projective flat 
$R$-scheme 
$Y$ whose completion is 
$Y_{/X}=\foX$. 
By Theorem \ref{functor smooth}, the functor 
$\LD_{X^\dagger}$ is formally smooth, so the 
$W[[T]]$-algebra 
$R$ is formally smooth.
Hence 
$X$ is smoothable and lifts to characteristic zero.

\textit{Step (ii): $X$ is of type III and in minus-one-form.}
For an invertible 
$\O_{X}$-module 
$\shL$, let 
$I_\shL\subset R$ be the smallest ideal so that  
$\shL$ admits an extension over 
$\foX\otimes R/I_\shL$. Since the obstruction group
$H^2(X,\O_{X})$ is 1-dimensional, we can find  an element
$r_\shL\in R$ contained in the maximal ideal 
$\maxid_R\subset R$ with 
$I_\shL=r_\shL R$ (compare \cite{Deligne 1976}, Prop.\ 1.5).

\begin{claim}
\label{claim}
There is an ample 
$\O_{X}$-module 
$\shL$ with 
$r_\shL\not\in \maxid^2_R + (2,T)$.
\end{claim}

Assume this for a moment.  By Theorem \ref{functor smooth}, the local ring
$R$ is  regular. Choose 
$\shL$ as in the  Claim.
Then the sequence 
$(2,T,r_\shL)$ is part of a regular system of parameters.
Hence the   map 
$W[[T]]\ra R/r_\shL R$ admits a section.
Using Grothendieck's Existence Theorem, we obtain a log smooth deformation 
$Y\ra\Spec(W[[T]])$. Consequently, 
$X$ is smoothable and lifts to characteristic zero.
Note that the restriction of
$Y$ over $ k[[T]]$ has smooth total space, and that the restriction over
$W$ is a locally trivial deformation.

It remains to prove the Claim.
Since 
$X$ is in  
minus-one-form, Proposition \ref{minus one projective} ensures that 
$X$ is projective.
Seeking a contradiction, we assume that 
$r_\shL\in \maxid^2_R + (2,T)$ for all ample 
$\O_{X}$-modules
$\shL$. This means that 
$\Pic(Y)\ra\Pic(X)$ is surjective for all locally trivial 
first-order log smooth deformations 
$Y^\dagger\ra \Spec(k[\epsilon]^\dagger)$.
By Serre duality and 
$\omega_X=\O_{X}$, the Yoneda pairing
$$
H^1(X,\Omega^1_{X^\dagger /k^\dagger}) \times 
H^1(X,\shHom(\Omega^1_{X^\dagger/k^\dagger},\O_{X}))\lra
H^2(X,\O_{X})
$$
in nondegenerate. A first-order locally trivial log deformation 
$Y^\dagger$ corresponds to an element 
$\zeta^\dagger\in H^1(\shHom(\Omega^1_{X^\dagger/k^\dagger},\O_{X}))$. 
Its image 
$\zeta\in H^1(\shHom(\Omega^1_{X/k},\O_{X}))$ corresponds to the 
underlying  deformation 
$Y$. According to \cite{Nakamura 1975}, Lemma 6.5, the obstruction map 
$\partial$ in 
$$
\Pic(Y)   \lra \Pic(X)  \stackrel{\partial}{\lra} H^2(X,\O_{X})
$$
factors over 
$$
\Pic(X)  \stackrel{\dlog}{\lra} H^1(X,\Omega_X^1) 
\stackrel{\langle\cdot,\zeta\rangle}{\lra}
H^2(X,\O_{X}).
$$
Since 
$\Pic(X)$ extends over all such deformations
$Y$, we infer that the composition 
$$
\Pic(X) \stackrel{\dlog}{\lra} H^1(X,\Omega^1_X)  \lra
H^1(X,(\Omega^1_{X/k})^{\vee\vee})   \lra
H^1(X,\Omega^1_{X^\dagger/k^\dagger}) 
$$ 
must be zero.
The commutative diagram in (\ref{diagram}) gives a commutative diagram 
$$
\begin{CD}
H^0(D,\shF)     @>>>  H^1(X,(\Omega^1_{X/k})^{\vee\vee}) 
@>>> H^1(X,\Omega^1_{X^\dagger/k^\dagger})   
\\ @VVV @VVV     @VVV\\
H^0(\tilde{C},\O_{\tilde{C}})     @>>>   H^1(S,\Omega^1_{S/k})              
@>>>
H^1(S,\Omega^1_{S^\dagger/k}).
\end{CD}
$$
with $
\shF=\ker(\O_{\tilde{D}}\ra\O_T) = \coker ((\Omega^1_{X/k})^{\vee\vee}\ra
\Omega^1_{X^\dagger/k^\dagger})$.
Since the map
$\Pic(X)\ra H^1(X,\Omega^1_{X^\dagger/k^\dagger})$ is zero,  the image of  
$\Pic(X)\ra H^1(S,\Omega^1_{S/k})$ is contained in the image of 
$H^0(D,\shF)\ra H^1(S,\Omega^1_{S/k})$.

We shall derive a contradiction as follows.
Since $X$ is in minus-one-form, the invertible $\O_S$-module $\O_S(C)$ 
descends to an invertible $\O_X$-module $\shL$ as in the 
proof of Proposition \ref{minus one projective}
Hence the Hodge class of $C$ is of the form $\sum\lambda_jC_j$
for certain coefficients $\lambda_j\in k$ coming from $H^0(D,\shF)$.
The latter condition means 
$\lambda_j=\lambda_k$ if the intersection
$\varphi(C_j)\cap\varphi(C_k)$ is a double curve, and 
$\lambda_j+\lambda_k=\lambda_l$ if the intersection
$\varphi(C_j)\cap\varphi(C_k)\cap\varphi(C_l)$ is a triple point.
According to Proposition \ref{linear independent}, there are at least 
two components, say 
$S_0,S_1\subset S$, with  
$\lambda_j=1$ for all 
$C_j\subset S_0\cup S_1$. 
Using the  double curve and triple point conditions, 
we infer that  the multiplicities in
$\sum\lambda_j C_j$ near $S_0$ are as follows:
\begin{center}
\unitlength1em
\begin{picture}(0,4)
\put(-2,2){\line(-1,1){2}}
\put(-2,2){\line(-1,-1){2}}
\put(-2,2){\line(1,0){4}}
\put(2,2){\line(1,-1){2}}
\put(2,2){\line(1,1){2}}
\put(-6,1.8){$S_0$}
\put(-2.5,3){$1$}
\put(-2.5,0){$1$}
\put(-.2,1){$0$}
\put(2,3){$0$}
\put(2,0){$0$}
\end{picture}
\end{center}
Inductively, we infer that the multiplicities $\lambda_j$ on each component 
$S_i\neq S_0$ are cycles of the form
$1,0,0,1,0,0,\ldots 1,0,0$. On the other hand, 
$\lambda_j=1$ for all
$C_j\subset S_1$, contradiction.

\textit{Step (iii): $X$ is of type III.}
By Proposition \ref{put in form}, there is a sequence of type I and type II 
modifications 
putting 
$X$ into 
minus-one-form 
$X'$. According to case (ii), there is a smoothing of 
$X'$ with smooth total space and a locally trivial deformation to
 characteristic zero.
Applying the reverse sequence of modifications to the total space of 
these deformations, 
we easily obtain the
desired deformations of 
$X$.

\textit{Step (iv):  $X$ is ordinary of type II.}
The arguments are similar to the case of type III surfaces in minus one form, 
so I leave them as an
exercise.

\textit{Step (v): $X$ is supersingular.}
Set 
$A=k[[t]]$.  
According to Lemma \ref{supersingular deforms}, there is a flat 
projective morphism 
$Y\ra\Spec(A)$ with closed fiber 
$Y_\sigma=X$, so that the generic fiber 
$Y_\eta$ is a d-semistable classical Enriques surface. 
Choose a closed embedding 
$Y\subset \PP^n_A$ for some 
$n>0$. Let 
$U\subset\Hilb_{\PP^n_\ZZ}$ be the Hilbert scheme of all 
normal crossing Enriques surfaces 
contained in 
$\PP^n_\ZZ$. We have just constructed a curve 
$\Spec(A)\subset U$. By case (i),
the point 
$u\in U$ representing the generic fiber 
$Y_\eta$ admits  generizations   corresponding to  
smooth Enriques surfaces and surfaces in
characteristic zero. Consequently,
$X$ is smoothable and
lifts to characteristic zero.
\qed


\end{document}